\documentclass[11pt]{article}
\usepackage{latexsym,amsfonts,amssymb,amsmath,amsthm}
\usepackage{epsfig}
\usepackage{color,comment}

\begin{document}

\newtheorem{thm}{Theorem}[section]
\newtheorem{cor}{Corollary}[section]
\newtheorem{lem}{Lemma}[section]
\newtheorem{prop}{Proposition}[section]
\newtheorem{defn}{Definition}[section]
\newtheorem{rk}{Remark}[section]
\newtheorem{nota}{Notation}[section]
\newtheorem{Ex}{Example}[section]
\def\nm{\noalign{\medskip}}

\numberwithin{equation}{section}

\newcommand{\ds}{\displaystyle}
\newcommand{\pf}{\medskip \noindent {\sl Proof}. ~ }
\newcommand{\p}{\partial}
\renewcommand{\a}{\alpha}
\newcommand{\z}{\zeta}
\newcommand{\pd}[2]{\frac {\p #1}{\p #2}}
\newcommand{\norm}[1]{\left\| #1 \right \|}
\newcommand{\dbar}{\overline \p}
\newcommand{\eqnref}[1]{(\ref {#1})}
\newcommand{\na}{\nabla}
\newcommand{\Om}{\Omega}
\newcommand{\ep}{\epsilon}
\newcommand{\tmu}{\widetilde \epsilon}
\newcommand{\vep}{\varepsilon}
\newcommand{\tlambda}{\widetilde \lambda}
\newcommand{\tnu}{\widetilde \nu}
\newcommand{\vp}{\varphi}
\newcommand{\RR}{\mathbb{R}}
\newcommand{\CC}{\mathbb{C}}
\newcommand{\NN}{\mathbb{N}}
\renewcommand{\div}{\mbox{div}~}
\newcommand{\bu}{{\bf u}}
\newcommand{\la}{\langle}
\newcommand{\ra}{\rangle}
\newcommand{\Scal}{\mathcal{S}}
\newcommand{\Lcal}{\mathcal{L}}
\newcommand{\Kcal}{\mathcal{K}}
\newcommand{\Dcal}{\mathcal{D}}
\newcommand{\tScal}{\widetilde{\mathcal{S}}}
\newcommand{\tKcal}{\widetilde{\mathcal{K}}}
\newcommand{\Pcal}{\mathcal{P}}
\newcommand{\Qcal}{\mathcal{Q}}
\newcommand{\id}{\mbox{Id}}
\newcommand{\stint}{\int_{-T}^T{\int_0^1}}

\newcommand{\be}{\begin{equation}}
\newcommand{\ee}{\end{equation}}

\newcommand{\rd}{{\mathbb R^d}}
\newcommand{\rr}{{\mathbb R}}
\newcommand{\alert}[1]{\fbox{#1}}
\newcommand{\eqd}{\sim}
\def\R{{\mathbb R}}
\def\N{{\mathbb N}}
\def\Q{{\mathbb Q}}
\def\C{{\mathbb C}}
\def\ZZ{{\mathbb Z}}
\def\l{{\langle}}
\def\r{\rangle}
\def\t{\tau}
\def\k{\kappa}
\def\a{\alpha}
\def\la{\lambda}
\def\De{\Delta}
\def\de{\delta}
\def\ga{\gamma}
\def\Ga{\Gamma}
\def\ep{\varepsilon}
\def\eps{\varepsilon}
\def\si{\sigma}
\def\Re {{\rm Re}\,}
\def\Im {{\rm Im}\,}
\def\E{{\mathbb E}}
\def\P{{\mathbb P}}
\def\Z{{\mathbb Z}}
\def\D{{\mathbb D}}
\def\p{\partial}
\newcommand{\ceil}[1]{\lceil{#1}\rceil}

\title{Lattice and Continuum Models Analysis of the Aggregation Diffusion Cell Movement}

\author{Lianzhang Bao\thanks{School of Mathematics, Jilin University, Changchun, 130012, P. R. China, and Department of Mathematics and Statistics,
Auburn University,  AL 36849, U. S. A. (lzbao@jlu.edu.cn) partially support by CPSF--183816.}\,\,   and
Zhengfang Zhou \thanks{Department of Mathematics, Michigan State University, East Lansing, MI 48824, USA (zfzhou@math.msu.edu).}}


\date{}

\maketitle

\begin{abstract}
The process by which one may take a discrete model of a biophysical process and construct a continuous model based on it is of mathematical interest as well as being of practical use.
In this paper, we first study the singular limit of a class of reinforced random walks on a lattice for which a complete analysis of the existence and stability of solutions are possible. In the continuous scenario, we obtain the regularity estimate of this aggregation diffusion model. As a by-product, nonexistence of solution of the continuous model with pure aggregation initial data is proved. When the initial is purely in diffusion region, asymptotic behavior of the solution is obtained. In contrast to continuous model, bounded-ness of the lattice solution, asymptotic behavior of solution in diffusion region with monotone initial date and the interface behaviors of the aggregation, diffusion regions are obtained. Finally we discuss the asymptotic behaviors of the solution under more general initial data with non-flux when the lattice points $N\leq 4$.
\end{abstract}


\medskip


\section{Introduction}
Movement is a fundamental process for almost all biological organisms, ranging from the single cell level to the population level.  There are three categories of  mathematical model that researchers have developed to describe this biological phenomena \cite{TSPS}. The first one is discrete model where space is divided into a lattice of points with the variables defined only at the points and time changes in "jump" \cite{BZ}, \cite{HPO}, \cite{S3}, \cite{T}. The second one is the continuous model where all variables are considered to be defined at every point in space and time changes continuously. The last one is the hybrid which is a mixture of the previous two \cite{KA1}, \cite{KC}. Each of these models have advantages and disadvantages depending on the phenomenon under consideration, and on the length scale over which we wish to investigate the phenomenon.

\par Discrete models of biophysical processes are of use when we are interested in the behavior of individual cells, as well as their interactions with other cells and the medium which surrounds them \cite{BZ}, \cite{HPO}. Different derivations may lead to very different behavior of the movement. For example, when the discrete population model of \cite{HPO} is considered: one can show that cells can move toward the high-density region and flee from low-density. But in \cite{BZ}, cells will move toward low-density and flee otherwise. Usually, the cells are considered to be points which move on a lattice according to certain rules. These rules can be modified according to the states of neighboring points, such as the density of the neighboring points \cite{BZ} or the volume and adhesion of the neighboring points \cite{KA1}, \cite{KC}. Individual based models have found useful application to many physical systems and even simple rules of interaction  can give rise to remarkably complex behavior. In particular, individual-based models have found applications in ecology, pattern formation, wound healing, tumor growth and gastrulation and vasculogenesis in the early embryo, amongst many others.

\par Continuous models frequently involve the development of a reaction-diffusion equation \cite{BZ}. These are useful when the length scale over which we wish to investigate the phenomenon is much greater than the diameter of the individual elements composing it. These models have been found to be particularly useful in the study of pattern formation in nature, especially the phenomenon of "diffusion driven instability" \cite{KA1}, \cite{KC}, \cite{BZ}, \cite{HPO}. The application of hybrid models - where cells are models as discrete entities with their movements being influenced by continuous spatial fields - has also been found to be a useful approach \cite{DS}. Discrete and continuous approaches to the modeling of cell migration developed by Bao and Zhou \cite{BZ} via biased random walk will be under investigated in the present paper.

\par If we consider both a continuous and discrete models of the same phenomenon, we would expect the models to give rise to similar solutions properties at length scales where their ranges of applicability overlap. If the solutions' behaviors are totally different, we should track their differences and find reasons in the original modeling.
\par The rest of the paper is organized as follows. In section 2, we consider the discrete model of biological cell movement and shown how it can be used to obtain an expression for the continuous aggregation diffusion model. Section 3 considers the existence and nonexistence of weak solution of the continuous models. Section 4 focuses on the analysis of the discrete model. Section 5 discusses the asymptotic behaviors of special case when $N\leq4$ and some open problems.

\section{Discrete and Continuous Models}
Recall the biological cell movement model \cite{BZ}.
By the need for survival, mating or to overcome the hostile environment, the population have the tendency of aggregation when the population density is small and diffusion otherwise. For simplicity, here we consider one species living in a one-dimensional habitat without birth term. To derive the model we follow a biased random walk approach plus a diffusion approximation. First we discretize space in a regular manner. Let $h$ be the distance between two successive points of the mesh and let $u(x,t)$ be the population density that any individual of the population is at the point $x$ and at time $t$. By scaling, we can assume that $0\leq u(x,t)\leq 1$. During a time period $\tau$ an individual which at time $t$ and at the position $x$, can either:
\begin{description}
  \item[1.] move to the right of $x$ to the point $x+h$, with probability $R(x,t),$ or
  \item[2.] move to the left of $x$ to the point $x-h$, with probability $L(x,t)$ or
  \item[3.] stay at the position $x$, with probability $N(x,t).$
\end{description}
Assume that there are no other possibilities of movement we have
\begin{equation*}
 N(x,t) + R(x,t) + L(x,t) = 1.
\end{equation*}
We assume that
\begin{eqnarray*}
R(x,t) &=& K(u(x+h,t)),
\\
L(x,t) &=& K(u(x-h,t)).
\end{eqnarray*}
Here we let $ K(u(x,t))$ measures the probability of movement which depend on the population.
\begin{figure}
\center
  \includegraphics[width=13cm]{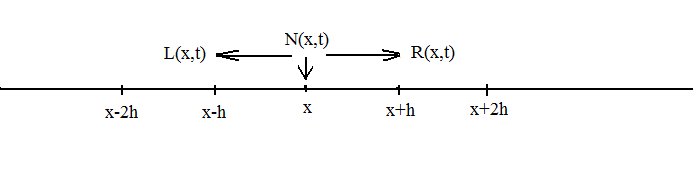}\\
  \caption{Movement of population}
\end{figure}
\par Using the notations above, the density $u(x,t)$ can be written as follows:
\begin{equation}\label{eq:001}
u(x,t+\tau) = N(x,t)u(x,t) + R(x-h,t)u(x-h,t) + L(x+h,t)u(x+h,t).
\end{equation}
By using Taylor series, we obtain the following approximation
\begin{eqnarray*}
u(x,t)+\tau\frac{du}{dt} &=&N(x,t)u(x,t)+\{R(x,t)u(x,t)-h\frac{d(R u)}{dx}+\frac{h^2}{2}\frac{d^2(R u)}{dx^2}\}
\\
&+&\{L(x,t)u(x,t) +h\frac{d(Lu)}{dx}+\frac{h^2}{2}\frac{d^2(Lu)}{dx^2}\},
\end{eqnarray*}
then we get
\begin{equation*}
\tau\frac{du}{dt} =\{-h\frac{d(Ru)}{dx}+\frac{h^2}{2}\frac{d^2(Ru)}{dx^2}\}+\{h\frac{d(Lu)}{dx}+\frac{h^2}{2}\frac{d^2(Lu)}{dx^2}\}.
\end{equation*}
Set
\begin{equation*}
\beta(x,t) = R(x,t)-L(x,t) = K(u(x+h,t))-K(u(x-h,t))= 2h\frac{d}{d x}[K(u(x,t)]+O(h^3)
\end{equation*}
and
\begin{equation*}
 \nu(x,t) =R(x,t)+L(x,t)= K(u(x+h,t))+K(u(x-h,t))= 2K(u(x,t))+O(h^2).
\end{equation*}
We can get
\begin{equation*}
 \tau\frac{d u}{d t}= -h\frac{d[(R-L)u]}{dx}+ \frac{h^2}{2}\frac{d^2[(R+L)u]}{d^2x}.
\end{equation*}
Now we substitute $\beta$ and $\nu$ in the above equation, we can get
\begin{equation*}
 \tau\frac{d u}{d t}= -2h^2\frac{d\{\frac{d}{d x}[K(u(x,t)]u\}}{dx}+ \frac{h^2}{2}\frac{d^2[2K(u(x,t))u]}{d^2x}+ O(h^3).
\end{equation*}
We assume that $h^2/\tau\rightarrow C>0$ (finite) as $\tau, h\rightarrow 0$, we get the following
\begin{equation*}
 \frac{d u}{d t}= -2C\frac{d\{\frac{d}{d x}[K(u(x,t)]u\}}{dx}+ C\frac{d^2[K(u(x,t))u]}{d^2x},
\end{equation*}
\begin{equation*}
 \frac{d u}{d t}= C\{-2\frac{d}{d x}[K(u(x,t)]u+ \frac{d[K(u(x,t))u]}{dx}\}_x,
\end{equation*}
\begin{equation*}
 \frac{d u}{d t}= C(-2K'u_xu+ K'u_xu+ Ku_x)_x= C[(K-uK')u_x]_x.
\end{equation*}
\par When higher order terms are kept, which will lead to the following
 \begin{equation}\label{eq:3001}
 u_t+\tau u_{tt}=((K-uK')u_x)_x+\frac{h^2}{12}[u_{xxxx}K-u\frac{d^4}{dx^4}(K(u))],
\end{equation}
the coefficient for the high order term $u_{xxxx}$ is $(K-uK')$ which is the same as the diffusion coefficient for $u_x$. Comparing to the standard Cahn-Hilliad equation or Cahn-Hilliad equation with degenerate mobility coefficient, our new equation \eqref{eq:3001} is totally new and the existence or nonexistence of the weak solution is still an open question.
\par If we let $K(u(x,t)) = 1/2(u^2-u^3)$, we can test that $0\leq K(u(x,t))\leq 1$ which satisfies the assumption of probability and that $K-uK'= u^2(u-1/2)$, which means aggregation when $0\leq u<1/2$. By plugging this probability in Equation \eqref{eq:001}, and we denote the discrete density $u(x,t)=u^t_j,u(x-h,t)=u^t_{j-1},u(x+h,t)=u^t_{j+1},\dots$ we have:
\begin{equation}\label{eq:002}
 u_j^{t+\tau}=u_j^t+\frac{u_j^tu_{j-1}^t}{2}(u_j^t+u_{j-1}^t-1)(u_{j-1}^t-u_j^t)+\frac{u_j^tu_{j+1}^t}{2}(u_j^t+u_{j+1}^t-1)(u_{j+1}^t-u_j^t),
\end{equation}
which is a special finite difference scheme of the following backward forward parabolic equation:
\begin{equation}\label{eq:003}
 u_t= [ D(u)u_x]_x \quad (x,t)\in Q_T,
\end{equation}
where $Q_T:=[0,1]\times[0,T]$,
\begin{equation}\label{eq:004}
 D(u)< 0 \quad \mbox{in} \quad(0,\alpha), \quad D(u)>0 \quad\mbox{in}\quad (\alpha,1).
\end{equation}

\section{Regularity and Properties of the Continuous Model}
In this section, we first give the definition of a weak solution of Equation \eqref{eq:003} with $D(u)$ always satisfies Equation \eqref{eq:004} and non-flux condition.
\begin{defn}\label{Def1}
 A locally continuous function $u(x,t)$ will be said to be a weak solution of the backward forward Equation \eqref{eq:003} with non-flux condition if
 \begin{equation}\label{eq:005}
 \int_0^1u^2(x,t)+|D(u)|u_x^2dx
 \end{equation}
 is uniformly bounded for $t\in [0,T]$ and if for any test function $\phi(x,t)$ in $C^1_0[Q_T]$,
 \begin{equation}\label{eq:006}
  \int\int[u\phi_t-D(u)u_x\phi_x]dxdt=0.
 \end{equation}
\end{defn}
We first establish the first derivative estimates. For simplicity, we limit ourselves to
showing how to obtain an a priori estimate. For any point $(x_0,t_0)$ such that $u(x_0,t_0)>\alpha$, there is a $\delta >0$ such that $D(u)>c>0$ for $ (x,t)\in \widetilde{Q}_\delta = [x_1,x_1+\delta]\times[t_0,t_0+\delta]$ and $(x_0,t_0)\in \widetilde{Q}_\delta$. We estimate
\begin{eqnarray*}
 \int\int\eta(x)^2u_{x}^2dxdt,
\end{eqnarray*}
where the cut-off function $\eta(x)$ is  smooth and vanishes
outside $[x_1,x_1+\delta]$. Differentiating the equation with respect to $t$ and
multiplying by $\eta^2u$ (if not differentiable, we use the smoothing operator to $u(x,t)$), we find
\begin{eqnarray*}
 \frac{1}{2}\frac{d}{dt}\int\eta^2u^2dx&=&\int\eta^2u_{t}udx=-\int
 D(u)u_x(\eta^2u)_x dx
\\
&=&\int-[D(u)\eta^2u_x^2+2D(u)\eta\eta_x u_x u]dx,
\end{eqnarray*}
So, integrating with respect to time and using Young inequality,
\begin{eqnarray*}
 c\int\int\eta^2u_{x}^2dxdt\leq
 C+\int\int2D(u)[\epsilon\eta^2u_{x}^2+(C/\epsilon)\eta_x^2u^2]dxdt,
\end{eqnarray*}
for some constant C and for every $\epsilon$. By taking $\epsilon$ small enough, we can get the first derivative estimate. Because $u$ is locally $H^1$ in $x$, $D(u)$ is strictly positive, by Nash-Moser estimate we can get $u\in C^{\alpha_1,\alpha_1/2},$ then according to Ladyzenskaya et al \cite{OVN}. [Chapter V, Theorem 7.4], we can get $u\in C^{1+\alpha_1,\alpha_1/2}.$

\par $C^{2+\alpha_2}$ estimate. Once we get the $C^{1+\alpha_1,\alpha_1/2}$ estimate, going back to our Equation \eqref{eq:003}, we have estimated that the equation for $u$ has H$\ddot{o}$lder continuous coefficients. If we view it as a linear equation, differentiating the equation, we now have
\begin{eqnarray*}
 u_{xt}&=&(D(u)u_{x})_{xx},
 \\
 &=& [D(u)u_{xx}+D'(u)u_x u_x]_x.
\end{eqnarray*}
If we let $v = u_x$, we get
\begin{eqnarray*}
 v_t=[D(u)v_x]_x+[D'(u)u_x]_xv+[D'(u)u_x]v_x,
\end{eqnarray*}
which is a uniformly forward parabolic equation for $v$. By using Schauder estimate, we can conclude that $u_x$ is locally $C^{2+\alpha_2,1+\alpha_2/2}$ for some $\alpha_2\in(0,1).$
\par Higher derivative estimates. Once $C^{2+\alpha_2,1+\alpha_2/2}$ estimates
have been found, it is easy to derive interior estimates for
higher-order derivatives. For instance, we may now differentiate the
equation once more to find an equation for $u_{xx}$:
\begin{eqnarray*}
 u_{xxt}-[D(u)u_{xx}+D'(u)u_x^2]_{xx}=0,
\end{eqnarray*}
which has
$H\ddot{o}lder$ continuous coefficients. Therefore, we conclude that
$u_{xx}$ is locally $C^{2+\alpha_2,1+\alpha_2/2}$, which gives a
fourth derivative estimate. The procedure is easily iterated and
continues to infinitum given $D(s)$ is infinitely differentiable.

\par From above estimate, we can get the following theorem for these big enough initial data.
\begin{thm}\label{thm1}
 Suppose that $D(u)$ satisfies \eqref{eq:004}. Then there exists a nonnegative classical solution in $C^{2,1}(Q_T)$ of Problem \eqref{eq:003} for all $T>0$, provided $\alpha<u(x,0)\in C^{1,\beta}([0,1]),\beta\in(0,1),$ and satisfies the non-flux conditions $D(u)u_x(0,t)=D(u)u_x(1,t)=0$.
\end{thm}
\par Furthermore, when $u_0\geq \alpha$ a.e. on $[0,1]$, we have the following asymptotic theorem
\begin{thm}\label{thm2}
 Suppose $u_0\geq \alpha, u_0\in C^{1,\beta}[0,1]$, and $u(x,t)$ is the solution of Equation \eqref{eq:003} with non-flux boundary condition, then solution $u(x,t)$ will go to constant $C=\int_0^1u(x,0)dx$.
\end{thm}
\proof First from the above estimate, we can see $u(x,t)$ is infinitely differentiable given $u(x,0)>\alpha$ a.e. on $[0,1]$. So we have
\begin{equation}\label{eq:007}
 \frac{d}{dt}\int_0^1u dx=\int_0^1u_t dx=\int_0^1(D(u)u_x)_xdx=D(u)u_x|_0^1=0,
\end{equation}
which gives the result of conservation of the total density.

\par When $u(x,0)>\alpha$ a.e. in $[0,1]$, we can have $D(u)>0$ by using the maximum principle of parabolic equation. Even more we can have $u(x,t)\geq \alpha+\delta_1$ for some small number $\delta_1$ and $t>0$.
\begin{equation}\label{eq:008}
 \frac{d}{dt}\int_0^1u^2 dx=2\int_0^1uu_t dx=2\int_0^1u(D(u)u_x)_xdx=-2\int_0^1D(u)u_x^2dx\leq 0,
\end{equation}
which means the "energy" of the solution is decreasing.
\par The last step we prove the solution will attend the constant $C$ as time $t$ goes to infinity.
\begin{eqnarray*}
 \frac{d}{dt}\int_0^1(u-C)^2dx &=&2\int_0^1(u-C)u_t dx =2\int_0^1(u-C)[D(u)u_x]_x dx
 \\
  &=&-2\int_0^1(u-C)_x D(u)u_xdx\leq -2\delta\int_0^1(u-C)_x^2dx
  \\
  &\leq&-2\delta c\int_0^1(u-C)^2dx \quad(\mbox{By using Poincar}\acute{e} \quad\mbox{inequality}).
\end{eqnarray*}
Then by using Gronwall's inequality, we can get
\begin{equation}\label{eq:009}
 \int_0^1(u-C)^2dx\leq (\int_0^1(u_0-C)^2dx)e^{-2\delta ct}.
\end{equation}
Let $t\rightarrow +\infty,$ we can get $\int_0^1(u-C)^2dx=0$ a.e., by combining $u$ is locally continuous, we obtain $\lim_{t\rightarrow+\infty}u(x,t)=C$.

\par The proof of $u-\alpha\geq\delta>0$ for all $t>0$. We claim
\begin{equation*}
\min_{x\in[0,1]}u(x,t)\geq \min_{x\in[0,1]}u_0(x),
\end{equation*}
otherwise
\begin{equation}\label{eq:010}
 \min_{x\in[0,1],0\leq t\leq T}u(x,t)=u(x_0,t_0),\quad T\geq t_0>0.
\end{equation}
Let $\nu=u(x,t)+\epsilon t,\epsilon>0$, then
\begin{eqnarray}\label{eq:011}
 \nu_t=u_t+\epsilon &=&[D(u)u_x]_x+\epsilon,\nonumber
 \\
 &=&[D(u)\nu_x]_x+\epsilon,\nonumber
 \\
 &=& D'(u)\nu_x^2+D(u)\nu_{xx}+\epsilon.
\end{eqnarray}
$\nu_x(0,t)=\nu_x(1,t)=0$, $\nu(x,t)$ can not choose a min for $t>0$. If not, $\nu(x_1,t_0)=\min.$
\begin{itemize}
  \item Case 1, $x_1\in(0,1),\nu_x=0,\nu_{xx}\geq 0,$ contradiction.
  \item Case 2, $x_1=0$ or 1, $\nu_x(0)=0,\nu_{xx}\geq 0$ contradiction.
\end{itemize}
So we obtain
\begin{equation*}
\min_{x\in[0,1],0\leq t\leq T}\nu(x,t)=\min_{x\in[0,1]}\nu(x,0)=\min_{x\in[0,1]}u(x,0),
\end{equation*}
 which means
\begin{equation}\label{eq:012}
 u(x,t)+\epsilon t\geq \min_{x\in[0,1]}u(x,0),\quad \forall \epsilon>0.
\end{equation}
When $\epsilon\rightarrow 0,u(x,t)\geq\min_{x\in[0,1]}u(x,0).$
\par For the case $u(x,0)\geq \alpha$, we use the approximation to Equation \eqref{eq:003}
\begin{equation}\label{eq:013}
 \left\{
   \begin{array}{ll}
     u_t^\epsilon=[(D(u)+\epsilon)u_x^\epsilon]_x, & (x,t)\in[0,1]\times (0,T] \\
     u_x^\epsilon(0,t)=u_x^\epsilon(1,t)=0, &  \\
     u^\epsilon(x,0)=u_0(x,0). &
   \end{array}
 \right.
\end{equation}
Using the above result, we can get $\lim_{t\rightarrow+\infty}u^\epsilon(x,t)=C$. Then let $\epsilon\rightarrow 0$, we can get $\lim_{t\rightarrow+\infty}u(x,t)=C$, for $u(x,0)\geq\alpha$ on [0,1].
\par\rightline{$\Box$}
By using the $C^\infty$ estimate above, we also can get the following non-existence result.
\begin{thm}\label{thm3}
 Suppose $0<u(x,0)<\alpha$ and $u(x,0)\in C^{1,\beta}[0,1]$, then Equation \eqref{eq:003} has no weak solution.
\end{thm}
\par
\proof If not, supposing there exists an weak solution $u(x,t)<\alpha$ in $Q_T$, by the linear transform $\tau=T-t$, we can reverse the time interval and get
\begin{equation}\label{eq:014}
 u_\tau=-(D(u)u_x)_x,\quad \forall (x,t)\in Q_T.
\end{equation}
Because $-D(u)>0$, then from the $C^\infty$ estimate, we can have $u(x,0)\in C^\infty(Q_T).$ If the initial value $u(x,0)$ is not infinitely differentiable, then there is no weak solution at all.
\par\rightline{$\Box$}

\section{Numerical analysis of the lattice model}
In our derivation of the new population model, we find every assumption is physically reasonable. But we find in Theorem \ref{thm3} that there is no weak solution when the initial date is not $C^\infty$ which lead us to consider the properties of the original lattice model and to see whether the original model is physically reasonable.
\begin{thm}\label{thm4}
 When considering equation \eqref{eq:002} with initial solution $1/2=\alpha\leq u(j,0)\leq 1,$ $ j=0,\dots,N$ and non-flux boundary condition, we obtain $1/2=\alpha\leq u(j,t)\leq 1$ for all $t\geq 0$. Furthermore, if $u(j,0)$ is monotone in $j$, then $u(j,t)$ is monotone in $j$ and
\begin{equation*}
  \lim_{t\rightarrow \infty} u(j,t)=\frac{1}{N-1}\sum_{j=1}^{N-1} u(j,0).
 \end{equation*}
\end{thm}
\begin{figure}
\center
  \includegraphics[width=12cm]{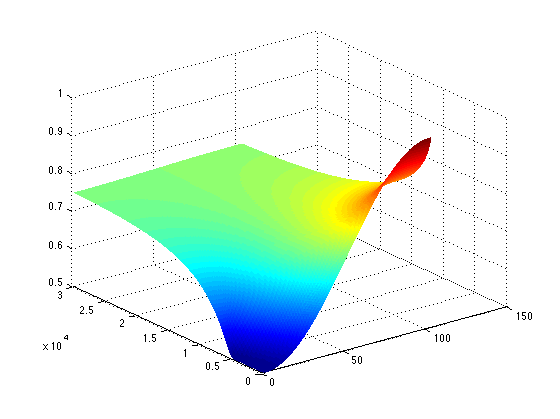}
  \caption{Asymptotic behavior of monotone initial solution}
\end{figure}
\proof Suppose $u(k_1,t)=\max_{0\leq j\leq N}u(j,t)$, $u(k_2,t)=\min_{0\leq j\leq N}u(j,t)$, then from Equation \eqref{eq:002} we can get
\begin{eqnarray}\label{eq:015}
 u(k_1,t+\tau)-u(k_1,t)&=&\frac{u_{k_1}^tu_{k_1-1}^t}{2}(u_{k_1}^t+u_{k_1-1}^t-1)(u_{k_1-1}^t-u_{k_1}^t)\nonumber
 \\
 &+&\frac{u_{k_1}^tu_{k_1+1}^t}{2}(u_{k_1}^t+u_{k_1+1}^t-1)(u_{k_1+1}^t-u_{k_1}^t),
\end{eqnarray}
$u_{k_1}^t+u_{k_1-1}^t-1\geq 0,u_{k_1}^t+u_{k_1+1}^t-1\geq 0$ and $u_{k_1-1}^t-u_{k_1}^t\leq0, u_{k_1+1}^t-u_{k_1}^t\leq 0$. Which lead to
\begin{equation*}
 u(k_1,t+\tau)-u(k_1,t)\leq 0.
\end{equation*}
The same idea can be used for the minimum value point and get
\begin{equation*}
 u(k_2,t+\tau)-u(k_2,t)\geq 0.
\end{equation*}
For arbitrary point $u(i,t)$
\begin{eqnarray*}
 u(i,t+\tau)-u(i,t)&=&\frac{u_i^tu_{i-1}^t}{2}(u_i^t+u_{i-1}^t-1)(u_{i-1}^t-u_i^t)+\frac{u_i^tu_{i+1}^t}{2}(u_i^t+u_{i+1}^t-1)(u_{i+1}^t-u_i^t)
\\
&\leq& 1/2[(u_{k_1}^t-u_i^t)]+1/2[(u_{k_1}^t-u_i^t)]
\\
&\leq& u_{k_1}^t-u_i^t.
\end{eqnarray*}
Again
\begin{eqnarray*}
 u(i,t+\tau)-u(i,t)&=&\frac{u_i^tu_{i-1}^t}{2}(u_i^t+u_{i-1}^t-1)(u_{i-1}^t-u_i^t)+\frac{u_i^tu_{i+1}^t}{2}(u_i^t+u_{i+1}^t-1)(u_{i+1}^t-u_i^t)
\\
&\geq& 1/2[(u_{k_2}^t-u_i^t)]+1/2[(u_{k_2}^t-u_i^t)]
\\
&\geq& u_{k_2}^t-u_i^t,
\end{eqnarray*}
which leads the result
\begin{equation*}
 \min_{0\leq j\leq N}u(j,0)\leq u(j,t)\leq\max_{0\leq j\leq N}u(j,0).
\end{equation*}

\par If the initial data is monotone and suppose
\begin{equation*}
 1/2=\alpha\leq u(0,0)= u(1,0)\leq u(2,0)\leq\dots\leq u(N-1,0)= u(N,0)\leq 1,
\end{equation*}
and let
\begin{equation}\label{eq:016}
 C_j^t=\frac{u_j^tu_{j-1}^t}{2}(u_j^t+u_{j-1}^t-1),
\end{equation}
we have
\begin{equation}\label{eq:017}
 \Delta U^{t+\tau}=[C]^t\Delta U^t,
\end{equation}
where
\begin{equation*}
 \Delta U^{t+\tau}=[(u_2^{t+\tau}-u_1^{t+\tau}),\dots,(u_N^{t+\tau}-u_{N-1}^{t+\tau})]^T,
\end{equation*}
\begin{equation*}
[C]^t=
\left(
         \begin{array}{ccccccccc}
           1-2C_1^t & C_2^t & 0 &  \dots & 0 & 0 & 0 \\
            C_1^t & 1-2C_2^t & C_3^t & \dots & 0 & 0 & 0  \\
           \vdots & \vdots & \vdots & \ddots & \vdots & \vdots & \vdots \\

           0 & 0 & 0 & \dots &C_{N-2}^t & 1-2C_{N-1}^t & C_{N}^t  \\
           0 & 0 & 0  & \dots & 0 &   C_{N-1}^t &1-2 C_N^t\\

         \end{array}
       \right)
\end{equation*}
When $1/2\leq u_j^t\leq 1$, $0\leq C_j^t=\frac{u_j^tu_{j-1}^t}{2}(u_j^t+u_{j-1}^t-1)\leq 1/2$, so $[C]^t$ is positive. By using the fact $U^t$ is positive, we obtain $U^{t+\tau}\geq0$. Initially $U^0\geq 0$, this lead to the conservation of monotonicity of $u_j^t$ in $0\leq j\leq N$ for all $t>0$.
\par From above two results, we can see
\begin{equation*}
\min_{0\leq j\leq N}u(j,t)=u(0,t)=u(1,t)
\end{equation*}
is bounded and increasing, so the limit
\begin{equation*}
 \lim_{t\rightarrow \infty}u(1,t+\tau)=\lim_{t\rightarrow \infty}u(1,t)+\lim_{t\rightarrow \infty}\frac{u_2^tu_1^t}{2}(u_1^t+u_2^t-1)(u_2^t-u_1^t)
\end{equation*}
exists. We can get $\lim_{t\rightarrow \infty}\frac{u_2^tu_1^t}{2}(u_1^t+u_2^t-1)(u_2^t-u_1^t)=0$ and $\frac{u_2^tu_1^t}{2}(u_1^t+u_2^t-1)>0$, which lead to $\lim_{t\rightarrow \infty}(u_2^t-u_1^t)=0$. Combined with the existence of $\lim_{t\rightarrow\infty}u_1^t$, we can get the existence of $\lim_{t\rightarrow\infty}u_2^t$. By using the same idea, which will lead to \begin{equation*}
 \lim_{t\rightarrow\infty}u_1^t=\lim_{t\rightarrow\infty}u_2^t=\dots=\lim_{t\rightarrow\infty}u_N^t=\frac{1}{N-1}\sum_{j=1}^{N-1}u(j,0).
\end{equation*}
\par\rightline{$\Box$}
\begin{thm}\label{thm5}
 For all the initial solution $0\leq u(j,0)\leq 1, j=0,\dots,N$, with non-flux boundary condition, the solution of Equation \eqref{eq:002} is bounded in $[0,1]$ and the total density is conservative.
\end{thm}
\proof For the non-flux boundary condition, we just assume
\begin{equation*}
 u(0,t)=u(1,t),\quad u(N,t)=u(N-1,t).
\end{equation*}
 From Equation \eqref{eq:002},we have
\begin{equation}\label{eq:018}
 u_j^{t+\tau}=u_j^t+C_j^t(u_{j-1}^t-u_j^t)+C_{j+1}^t(u_{j+1}^t-u_j^t).
\end{equation}
We get
\begin{eqnarray*}
 \sum_{j=1}^{N-1} u_j^{t+\tau}=\sum_{j=1}^{N-1} u_j^{t}&+&C_1^t(u_0^t-u_1^t)+C_2^t(u_2^t-u_1^t)
 \\
 &+&C_2^t(u_1^t-u_2^t)+C_3^t(u_3^t-u_2^t)+\dots+
 \\
 &+&C_{N-1}^t(u_{N-2}^t-u_{N-1}^t)+C_N^t(u_N^t-u_{N-1}^t),
\end{eqnarray*}
the summation of the $C_j^t$ terms will result to zero. From which we obtain
\begin{equation}\label{eq:017}
 \sum_{j=1}^{N-1} u_j^{t+\tau}=\sum_{j=1}^{N-1} u_j^{t}
\end{equation}
and the lattice model is conservative.
\par When rewriting Equation \eqref{eq:002} to the form
\begin{equation*}
 u_j^{t+\tau}=\frac{u_j^t}{2}+\frac{u_j^tu_{j-1}^t}{2}(u_j^t+u_{j-1}^t-1)(u_{j-1}^t-u_j^t)+\frac{u_j^t}{2}+\frac{u_j^tu_{j+1}^t}{2}(u_j^t+u_{j+1}^t-1)(u_{j+1}^t-u_j^t),
\end{equation*}
we consider the following equation
\begin{equation}\label{eq:018}
 f(x,y)=x+xy(x+y-1)(y-x), \quad x,y\in [0,1].
\end{equation}
If we can find the upper and lower bound of $f(x,y)$ such that $0\leq f(x,y)\leq 1$, then we obtain the result $u(j,t)\in [0,1]$ for all $t>0$.
\begin{eqnarray*}
 &&f(0,y)=0,
 \\
 &&f(1,y)=1+y^3-y^2, f_y(1,y)=3y^2-2y.
\end{eqnarray*}
When we let $f_y(1,y)=0, y=2/3$ or 0 which lead to $f(1,2/3)=23/27>1/2, f(1,0)=1$ and $1/2<f(1,y)\leq 1$. Again, we test other boundary value $f(x,0)=x$ and $f(x,1)=x+x^2(1-x),f_x(x,1)=1+2x-3x^2$. Let $f_x(x,1)=1+2x-3x^2=0$, we get $x=-1/3$ or $x=1$ and $0\leq f(x,1)\leq 1$. Now we test if the function can get the extreme value when $x,y\in(0,1)$. Taking the derivative, we get
\begin{eqnarray*}
 f_x&=& 1+y(x-y+y^2-x^2)+xy(1-2x)=0,
 \\
 f_y&=&x(x-y+y^2-x^2-y+2y^2)=0.
\end{eqnarray*}
After simplification and let $x\neq 0$, we get
\begin{eqnarray*}
 f_x&=& 1+y^2(1-2y)+xy(1-2x)=0,
 \\
 f_y&=&3y^2-x^2+x-2y=0.
\end{eqnarray*}
We can obtain $f_x>0$ for $x\in(0,1/2)$ or $y\in(0,1/2)$, when $x,y\in[1/2,1]$, by using Theorem \ref{thm4}, we have $1/2\leq u(j,t)\leq 1$ for any $0\leq j\leq N, t\geq 0$ when $0\leq u(j,0)\leq 1$. But from equation $f_x=0$ which mean the maximum or minimum points, we have
\begin{equation*}
 0\leq u(j,t)\leq 1, \quad \forall 0\leq j\leq N, t\geq 0.
\end{equation*}
\par\rightline{$\Box$}
In the following we consider the interface behaviors of the forward region
\begin{equation*}
Q_d^+(t)=\{(i,t)| u(i,t)\geq 1/2,i\in [1,2,3,\dots,N]\},
\end{equation*}
backward region
\begin{equation*}
Q_d^-(t)=\{(i,t)| u(i,t)< 1/2,i\in [1,2,3,\dots,N]\}.
\end{equation*}
\begin{thm}\label{thm6}
 Under any general initial condition $0\leq u(j,0)\leq 1$, for any $t_1>t\geq 0$
\begin{equation}\label{eq:019}
 Q_d^+(t)\subseteq Q_d^+(t_1).
\end{equation}
\end{thm}
\proof We first begin at the simplest situation when $Q_d^+(t)=\{(i,t)| u(i,t)\geq 1/2,i\in [1,2,3,\dots,N]\}=(k,t)$. If $u(k,t)=1/2$, from Equation \eqref{eq:002}
\begin{equation*}
 u_k^{t+\tau}=u_k^t+\frac{u_k^tu_{k-1}^t}{2}(u_k^t+u_{k-1}^t-1)(u_{k-1}^t-u_k^t)+\frac{u_k^tu_{k+1}^t}{2}(u_k^t+u_{k+1}^t-1)(u_{k+1}^t-u_k^t),
\end{equation*}
and we have $(u_k^t+u_{k-1}^t-1)<0, (u_{k-1}^t-u_k^t)<0, (u_k^t+u_{k+1}^t-1)<0,(u_{k+1}^t-u_k^t)<0$, which lead to
\begin{equation*}
 u_k^{t+\tau}\geq u_k^t= 1/2.
\end{equation*}
Suppose $u(k,t)>1/2$, then we have $u(k,t)=1/2+\delta$ for some positive $\delta$ and Equation \eqref{eq:002} is of the form
\begin{equation*}
 u_k^{t+\tau}=1/2+\delta+\frac{u_k^tu_{k-1}^t}{2}(u_k^t+u_{k-1}^t-1)(u_{k-1}^t-u_k^t)+\frac{u_k^tu_{k+1}^t}{2}(u_k^t+u_{k+1}^t-1)(u_{k+1}^t-u_k^t).
\end{equation*}
If $(u_k^t+u_{k-1}^t-1)\leq0$ and $(u_k^t+u_{k+1}^t-1)\leq0$, it is easy to obtain $u(k,t+\tau)\geq u(k,t)$. Suppose $(u_k^t+u_{k-1}^t-1)>0$ and $(u_k^t+u_{k+1}^t-1)<0$, we have $u(k-1,t)>1-u(k,t)=1/2-\delta$ and
\begin{equation*}
 u(k,t)+u(k-1,t)-1<1/2+\delta+1/2-1<\delta,\quad u(k,t)-u(k-1,t)<2\delta.
\end{equation*}
Then
\begin{equation*}
 u_k^{t+\tau}>1/2+\delta-\delta\{u_k^tu_{k-1}^t\}\delta+\frac{u_k^tu_{k+1}^t}{2}(u_k^t+u_{k+1}^t-1)(u_{k+1}^t-u_k^t).
\end{equation*}
Because $u_k^tu_{k-1}^t<1/4$ and $0<\delta<1/2$, we can obtain
\begin{equation}\label{eq:021}
 \delta- 2\delta^2u_k^tu_{k-1}^t\geq 0,
\end{equation}
which lead to $u(k,t+\tau)\geq u(k,t)$ even for both $(u_k^t+u_{k-1}^t-1)>0,(u_k^t+u_{k+1}^t-1)>0$.
\par When $Q_d^+(t)=\{(i,t)| u(i,t)\geq 1/2,i\in [1,2,3,\dots,N]\}$ has more than two successive points, we rewrite Equation \eqref{eq:002} as
\begin{equation}\label{eq:022}
 u_j^{t+\tau}=\frac{u_j^t}{2}+\frac{u_j^tu_{j-1}^t}{2}(u_j^t+u_{j-1}^t-1)(u_{j-1}^t-u_j^t)+\frac{u_j^t}{2}+\frac{u_j^tu_{j+1}^t}{2}(u_j^t+u_{j+1}^t-1)(u_{j+1}^t-u_j^t).
\end{equation}
Let's consider the boundary points with $0\leq u(j-1,t)<1/2$ and $1/2\leq u(j,t),u(j+1,t)\leq 1$. Using the result in simplest situation
\begin{equation}\label{eq:023}
 \frac{u_j^t}{2}+\frac{u_j^tu_{j-1}^t}{2}(u_j^t+u_{j-1}^t-1)(u_{j-1}^t-u_j^t)\geq 1/2\times1/2=1/4,
\end{equation}
and if $ u(j,t)\geq u(j+1,t)\geq 1/2$
\begin{eqnarray*}
 &&\frac{u_j^t}{2}+\frac{u_j^tu_{j+1}^t}{2}(u_j^t+u_{j+1}^t-1)(u_{j+1}^t-u_j^t)
 \\
 &&> 1/2[1/2+\delta+ (1+1-1)(1/2-1/2-\delta)]=1/2(1/2+\delta-\delta)=1/4,
\end{eqnarray*}
otherwise if $ u(j,t)\leq u(j+1,t)$,
\begin{eqnarray*}
 &&\frac{u_j^t}{2}+\frac{u_j^tu_{j+1}^t}{2}(u_j^t+u_{j+1}^t-1)(u_{j+1}^t-u_j^t)
 \\
 &&\geq \frac{1}{2}u_j^t\geq 1/4.
\end{eqnarray*}
When considering the points inside the forward region, we can have the result that $1/2\leq u(j,t)\leq 1$ by using the maximum and minimum principle from Theorem \ref{thm4}. Combining all these together, we complete the proof of theorem.
\par\rightline{$\Box$}
\section{Asymptotic behaviors under special case when $N\leq4$}
\par In this section, we consider a special case when $N=4$ and $u(0,t)=u(4,t)=0$ for $t\geq 0$ which corresponding to the hostile environment in biology, when $N=1,2,3$, the asymptotic behaviors of the solutions are easy to be obtained.
\par Case 1: $u(1,0)+u(2,0)>1, u(2,0)+u(3,0)<1$ and $0\leq u(1,0)<u(2,0)$. From Equation \eqref{eq:002}, we have
\begin{eqnarray}
 u_1^{t+\tau}&=&u_1^t+\frac{u_1^tu_{2}^t}{2}(u_1^t+u_{2}^t-1)(u_{2}^t-u_1^t),\label{eq:3025}
 \\
 u_2^{t+\tau}&=&u_2^t+C_1^t(u_{1}^t-u_2^t)+C_2^t(u_{3}^t-u_2^t),\label{eq:3026}
 \\
 u_1^{t+\tau}+u_2^{t+\tau}&=&u_1^t+u_2^t+\frac{u_2^tu_{3}^t}{2}(u_2^t+u_{3}^t-1)(u_{3}^t-u_2^t),\label{eq:3027}
 \\
 u_2^{t+\tau}-u_1^{t+\tau}&=&(u_2^t-u_1^t)(1-2C_1^t)+\frac{u_2^tu_{3}^t}{2}(u_2^t+u_{3}^t-1)(u_{3}^t-u_2^t),\label{eq:3028}
 \\
 u_3^{t+\tau}&=&u_3^t+\frac{u_2^tu_{3}^t}{2}(u_2^t+u_{3}^t-1)(u_{2}^t-u_3^t),\label{eq:3029}
 \\
 u_2^{t+\tau}+u_3^{t+\tau}&=&u_2^t+u_3^t+\frac{u_2^tu_{1}^t}{2}(u_2^t+u_{1}^t-1)(u_{1}^t-u_2^t).\label{eq:3030}
\end{eqnarray}
If the initial solution satisfy case 1, we can see $u_1^t,u_1^{t}+u_2^{t}$ are increasing, $u_3^t,u_2^{t}+u_3^{t}$ are decreasing. Also from Theorem \ref{thm4}, $u(j,t)$ are bounded in $[0,1]$ which lead to the existence of the limit of each term and the following asymptotic behaviors of the solution
\begin{eqnarray*}
 \lim_{t\rightarrow \infty}u_1^{t}&\rightarrow&\overline{u}_1=\lim_{t\rightarrow \infty}u_2^{t}\rightarrow\overline{u}_2>1/2,
 \\
 \lim_{t\rightarrow \infty}u_3^{t}&\rightarrow&\overline{u}_3=0.
\end{eqnarray*}

\par Case 2: $u(1,0)+u(2,0)<1, u(2,0)+u(3,0)<1$ and $u(2,0)$ is the min. From Equations \eqref{eq:3025}-\eqref{eq:3030}, we have
 $u_1^t,u_3^{t}$ are increasing, $u_2^t,u_1^{t}+u_2^{t},u_2^{t}+u_3^{t}$ are decreasing and $u(j,t)$ are bounded in $[0,1]$ which lead to the existence of the limit of each term and the following asymptotic behaviors of the solution
\begin{eqnarray*}
 \lim_{t\rightarrow \infty}u_1^{t}&\rightarrow&\overline{u}_1,
 \\
 \lim_{t\rightarrow \infty}u_2^{t}&\rightarrow&\overline{u}_2=0,
 \\
 \lim_{t\rightarrow \infty}u_3^{t}&\rightarrow&\overline{u}_3.
\end{eqnarray*}

\par Case 3: $u(1,0)+u(2,0)<1, u(2,0)+u(3,0)<1, u(1,0)+u(2,0)+u(3,0)<1$ and $u(2,0)$ is the maximum. From Equations \eqref{eq:3025}-\eqref{eq:3030}, we have
 $u_1^t,u_3^{t}$ are decreasing, $u_2^t$ is increasing and $u(j,t)$ are bounded in $[0,1]$ which lead to the existence of the limit of each term and the following asymptotic behaviors of the solution
\begin{eqnarray*}
 \lim_{t\rightarrow \infty}u_1^{t}&\rightarrow&\overline{u}_1=0,
 \\
 \lim_{t\rightarrow \infty}u_2^{t}&\rightarrow&\overline{u}_2=u(1,0)+u(2,0)+u(3,0),
 \\
 \lim_{t\rightarrow \infty}u_3^{t}&\rightarrow&\overline{u}_3=0.
\end{eqnarray*}
\par Case 4: $u(1,0)+u(2,0)<1, u(2,0)+u(3,0)<1, u(1,0)+u(2,0)+u(3,0)<1$ and $u(3,0)$ is the maximum, $u(1,0)$ is the minimum, otherwise it is the case 2. (When $u(1,0)$ is the maximum, it the same as $u(3,0)$ is the maximum.) From Equations \eqref{eq:3025}-\eqref{eq:3030}, we assume that $u^{t+\tau}_3>u^{t+\tau}_2>u^{t+\tau}_1$, otherwise it goes to case 2 or case 3, the limits will exist. In this case $u_3^t$ is increasing and $u(j,t)$ are bounded in $[0,1]$ which lead to the existence of the limit of each term and the following asymptotic behaviors of the solution
\begin{eqnarray*}
 \lim_{t\rightarrow \infty}u_1^{t}&\rightarrow&\overline{u}_1=0,
 \\
 \lim_{t\rightarrow \infty}u_2^{t}&\rightarrow&\overline{u}_2=0,
 \\
 \lim_{t\rightarrow \infty}u_3^{t}&\rightarrow&\overline{u}_3=u(1,0)+u(2,0)+u(3,0).
\end{eqnarray*}
\par Case 5: $u(1,0)+u(2,0)<1, u(2,0)+u(3,0)<1,1\leq u(1,0)+u(2,0)+u(3,0)$ and $u(3,0)$ is the maximum. (When $u(1,0)$ is the maximum, it the same as $u(3,0)$ is the maximum.)
\par If $u(2,0)$ is the minimum, it goes to case 2. So we consider the case when $u(3,0)>u(2,0)>u(1,0)$. From Equations \eqref{eq:3025}-\eqref{eq:3030}, we have
 $u_1^t,u_1^{t}+u_2^{t}$ are decreasing, $u_3^t,u_2^{t}+u_3^{t}$ are increasing first. Because $1\leq u(1,0)+u(2,0)+u(3,0)$, when $ u_2^{t}+u_3^{t}$ first pass 1 and $u(2,t)$ become maximum, then this become similar to Case 1 and the solution $u(j,t)$ have the asymptotic limits. When $u(1,t)+u(2,t)<1, u(2,t)+u(3,t)>1$ and $u(3,t)$ is the maximum, in this case $u(1,t)+u(2,t)<1$ for all $t$, $u(1,t)$ is decreasing and $u(2,t)>u(1,t)$ for all $t$ and we can get
\begin{equation*}
\lim_{t\rightarrow \infty}u_1^t=0, \lim_{t\rightarrow\infty}u_2^t=\lim_{t\rightarrow}u_3^t=1/2\{u(1,0)+u(2,0)+u(3,0)\}.
\end{equation*}

\par Case 6: $u(1,0)=u(3,0)$ and $u(1,0)+u(2,0)>1$. If $u(2,0)$ is the minimum, then $u(2,t), u(1,t)+u(2,t), u(2,t)+u(3,t)$ are increasing and $u(1,t), u(3,t)$ are decreasing. From $\lim_{t\rightarrow\infty}u(1,t)+u(2,t)>1$ exists and $u(2,t), u(1,t)>0$ for all $t$, we can get the asymptotic behaviors of solutions
 \begin{equation*}
\lim_{t\rightarrow\infty}u(1,t)=\lim_{t\rightarrow\infty}u(2,t)=\lim_{t\rightarrow\infty}u(3,t)=1/3(u(1,0)+u(2,0)+u(3,0)).
\end{equation*}
When $u(2,0)$ is the maximum point, in this case $u(1,t), u(3,t)$ are increasing and $u(2,t),u(1,t)+u(2,t)$ are decreasing. We suppose that $u_1^t+u_2^t=1+\delta,$ then
\begin{equation*}
 u_1^{t+\tau}+u_2^{t+\tau}=1+\delta - \delta\frac{u_1^tu_2^t}{2}(u_2^t-u_1^t)>1,
\end{equation*}
$u(1,t)=u(3,t)>1/2$ for $t>t_0$. In this case, it will go to special case 10 and we obtain
\begin{equation*}
 \lim_{t\rightarrow\infty}u_1^t=\lim_{t\rightarrow\infty}u_2^t=\lim_{t\rightarrow\infty}u_3^t.
\end{equation*}

\par Case 7: $u(1,0)+u(2,0)>1$ and $u(2,0)+u(3,0)>1$. If $u(2,0)$ is the minimum point, in this case, we can get $u_1^t+u_2^t>1, u_2^t+u_3^t>1$ for all $t>0$. So we can get the following
\begin{eqnarray*}
 |u_2^{t+\tau}-u_1^{t+\tau}|&\leq& \max\{|1-2C_1^t||u_2^t-u_1^t|, C_2^t|u_3^t-u_2^t|\}
 \\
  |u_3^{t+\tau}-u_2^{t+\tau}|&\leq& \max\{|1-2C_2^t||u_3^t-u_2^t|, C_1^t|u_2^t-u_1^t|\},
\end{eqnarray*}
which leads to the convergence of the asymptotic behaviors of $u(1,t),u(2,t)$ and $u(3,t)$.
\par In the case $u(3,0)>u(2,0)>1/2>u(1,0)$, we can get $u(1,t)+u(2,t)$ is increasing first and $u(2,t)+u(3,t)$ is decreasing but always larger than 1. From theorem \ref{thm4}, if $u(3,t)$ keep as maximum in this situation, then we can see the existence of the limits of the solution. But when $u(2,t), u(3,t)$ change their orders, we can find $u(1,t)+u(2,t)$ is decreasing and $u(1,t)+u(2,t)$ may be less than 1, in this case we still have some difficulties to prove if the system have asymptotic limits or periodic solutions. This is one of the open questions we still need to answer.

\par Case 8: $u(1,0)+u(2,0)<1$ and $u(2,0)+u(3,0)>1, u(2,0)>1/2$ and $u(3,0)$ is the maximum. We can see $u(1,t)+u(2,t)$ is increasing and $u(3,t)$ is decreasing. If $u(1,t)+u(2,t)<1$ for all $t>0$, then we can get the limits of the solution. But if $u(1,t_0)+u(2,t_0)>1$ and $u(2,t_0)$ is the maximum for some $t_0>0$, this goes back to the open question we mentioned above.
\par Case 9: $u(1,0)+u(2,0)=1, u(2,0)=u(3,0)$, by using Equations \eqref{eq:3025}-\eqref{eq:3030}, we can see that this is an trivial steady state solution.
\par Case 10: $1/2\leq u(1,0),u(3,0)\leq u(2,0)\leq 1$ or $1/2\leq u(2,0)\leq u(1,0),u(3,0)\leq 1$. For the first situation $u(2,0)$ is the local maximum. Using Equations \eqref{eq:3025}-\eqref{eq:3030}, we have
 $u_1^{t+\tau}\geq u_1^t,u_3^{t+\tau}\geq u_3^t,u_2^{t+\tau}\leq u_2^t$ and
 \begin{eqnarray*}
  u_2^{t+\tau}-u_1^{t+\tau}&=&(u_2^t-u_1^t)(1-2C_1^t)+\frac{u_2^tu_{3}^t}{2}(u_2^t+u_{3}^t-1)(u_{3}^t-u_2^t)
  \\
  &\leq &(u_2^t-u_1^t)(1-2C_1^t)<(u_2^t-u_1^t).
 \end{eqnarray*}
If in the case $u_2^{t+\tau}-u_1^{t+\tau}<0$, we can check $|u_2^{t+\tau}-u_1^{t+\tau}|<1/2|u_2^t-u_3^t|$.
But when $u(2,t)$ is the local minimum
  \begin{eqnarray*}
  u_1^{t+\tau}-u_2^{t+\tau}&=&(u_1^t-u_2^t)(1-2C_1^t)-\frac{u_2^tu_{3}^t}{2}(u_2^t+u_{3}^t-1)(u_{3}^t-u_2^t)
  \\
  &\leq &(u_1^t-u_2^t)(1-2C_1^t).
  \\
  u_2^{t+\tau}&\leq& u_2^t + 1/2(u_1^t-u_2^t) +1/2(u_3^t-u_2^t)\leq \max\{u_1^t,u_3^t\}.
 \end{eqnarray*}
Using the above three situations, we can obtain
\begin{equation}
 \max (|u_3^{t+\tau}-u_2^{t+\tau}|,|u_1^{t+\tau}-u_2^{t+\tau}|)\leq \max[(1-2C_1^t),1/2]\{ \max (|u_3^{t}-u_2^{t}|,|u_1^{t}-u_2^{t}|)\},
\end{equation}
where $C_1^t=\frac{u_1^tu_2^t}{2}(u_1^t+u_2^t-1)\leq \max\{u(1,0),u(2,0),u(3,0)\}.$ So $0<(1-2C_1^t)<1$ for $u(1,t)\neq u(2,t)$ (if $u(1,0)=u(2,0)$ and it is easy to test the sequence will change to monotone sequence and which will converge to the average), which lead to
 \begin{equation*}
  \lim_{t\rightarrow \infty}|(u_1^t-u_2^t)|=0,\quad \lim_{t\rightarrow\infty}|(u_3^t-u_2^t)|=0.
 \end{equation*}
We get
 \begin{equation*}
  \lim_{t\rightarrow \infty}u(1,t)=\lim_{t\rightarrow \infty}u(2,t)=\lim_{t\rightarrow \infty}u(3,t).
\end{equation*}
\par\rightline{$\Box$}
We discussed the asymptotic behaviors of the lattice model when the initial data is monotone in the diffusion region and in the case $N\leq4$. In our simulation (see FIG. 5.1), we find the asymptotic behaviors of the lattice model also converge to the average of the initial solution without monotonicity assumption and totally in diffusion region. Whether this lattice solution is periodic or converge to some constant is still under investigation.
%

\end{document}